\documentclass[12pt]{article}

\def\equ{\equiv \kern-.25cm \not}
\font\gordas = msbm10 at 12pt
\def\cajita{\rule{5pt}{5pt}}
\font\gorditas = msbm10 at 8pt
\def\bbb#1{\hbox {{\gordas #1}}}
\def\bbba#1{\hbox {{\gorditas #1}}}

\def\errita{{\bbba R}}
\def\erre{{\bbb R}}

\def\Dom{{\rm Dom}}

\newtheorem{theorem}{Theorem}
[section]
\newtheorem{proposition}[theorem]
{Proposition}
\newtheorem{lemma}[theorem]{Lemma}
\par 

\par 
\newtheorem{corolary}[theorem]{Corollary}
\par 

\par 

\par

\par

\par 

\par 

\par 

\par 

\begin{document}

\begin{center}
{\Large\bf An anticipating It\^o formula for L\'evy processes}\\[1cm]
{\bf Elisa Al\`{o}s$^*$}\\
Dpt. d'Economia i Empresa\\
Universitat Pompeu Fabra\\
c/Ramon Trias Fargas, 25-27\\
08005 Barcelona, Spain\\[.5cm]
{\bf Jorge A. Le\'{o}n$^\dagger$}\\
Depto. de Control Autom\'{a}tico\\
Cinvestav-IPN\\
Apartado Postal 14-740\\
07000 M\'{e}xico, D.F., Mexico\\[.5cm]
{\bf Josep Vives$^\ddagger$}\\
Dpt.de  Probabilitat, L\`{o}gica i Estad\'{\i}stica\\
 Universitat de Barcelona\\
Gran Via, 585\\
 08007 Barcelona, Spain
\end{center}
\footnote{\kern -.75cm $^*$Supported by grants CICYT\ BFM2003-04294 and MCYT SEC2003-04476\\
$^\dagger$Partially supported by the CONACyT grant 45684-F \\
$^\ddagger$Supported by grant MECD-BFM2003-00261} 

\begin{abstract}
In this paper, we use the Malliavin calculus techniques to obtain an 
anticipative
version of the change of variable formula for  L\'evy processes. Here 
the coefficients are in the domain of
the anihilation (gradient) operator in the ``future sense'',
which includes the family of all adapted and square-integrable
processes. This domain was introduced on the Wiener space by
Al\`os and Nualart \cite{AN}. Therefore, our It\^o formula is not
only an extension of the usual adapted formula for L\'evy processes, but
also an extension of the anticipative version on Wiener space obtained 
in \cite{AN}.

 \end{abstract}



\section{Introduction}

  It is well--known that the It\^o formula, or change of
 variable formula, is one of the most powerful tools of the stochastic analysis
 due to its vast  range of applications. So, in the  last few years, various
 researchers have studied extensions of the classical It\^o formula for
 different interpretations of stochastic integral 
(see, for instance, Al\`os and Nualart \cite{AN}, Di Nunno et al. 
\cite{Nupro}, Moret and Nualart \cite{MoN}, Nualart and Taqqu \cite {NT},
and Tudor and Viens \cite {TuV}).  In particular, several authors
have been interested in finding extensions of this important formula
to the case where the coefficients are not adapted to the underlying
filtration (see Di Nunno et al \cite{Nume}, Le\'on et al. \cite{JRD}, 
Nualart and Pardoux \cite{NP}, or Russo and Vallois \cite{RV}).

The Malliavin calculus or calculus of variations is another 
 important tool of the stochastic analysis that allows us to deal with 
stochastic integrals whose domains include processes that are not necessarily
 adapted to the underlying filtration. Recently, the interest of this calculus
 has increased considerably because of its applications in finance 
(see, for example, Al\`os \cite{alo}, Al\`os et al.  \cite{ALV}, 
Bally et al. \cite{BC},
Fourni\'e et al. \cite{Fou1,Fou2}, Imkeller \cite{Imk},
Nualart \cite{Nul} or {\O}ksendal \cite{Oks}), or other
 theoretical applications
 (see Al\`os and Nualart \cite{AN}, Le\'on and Nualart
\cite{LeN}, Nualart \cite{Nua,Nul} or Sanz-Sol\'e \cite{Sanz}).
 This important theory is basically based on the
 divergence and gradient operators.

The divergence operator has been interpreted as a stochastic integral because 
it has properties similar to those of the It\^o stochastic integral. For 
instance, the isometry and local properties, 
the fact that it can be approximated  by 
Riemann sums, the integration by parts formula, etc. (see Nualart \cite{Nul}).
Hence, it is important to count on a change of variable formula for
the divergence operator in order to improve
the applications of the Malliavin calculus to  different
areas of the human knowledge. 

On the Wiener space, the divergence operator was defined by Skorohod 
\cite{Sko} and it is an extension of the classical It\^o integral. In 
order to analyze the properties of the Skorohod integral, the adaptability of
 the integrands (necessary in the It\^o's calculus) is changed by some analytic
 properties that are used to define some spaces, called Sobolev spaces, 
where a fundamental ingredient is the derivative (gradient) operator
(see Sections \ref{sub:2.3} and \ref{sub:2.4} below). For 
instance, Al\'os and Nualart \cite{AN} have considered  processes with 
derivatives ``in the future sense''.  
In this paper the stochastic integral with respect to the continuous part
of the underlying L\'evy process is in the Skorohod sense.
The Skorohod integral can be introduced 
using  different approaches. Namely, the first method is via the Wiener chaos 
decomposition, and the second one considers the Skorohod integral as the 
adjoint of the gradient (derivative) operator.

On the Poisson space, the above two methods produce different definitions of 
stochastic integral (see, for example Carlen and Pardoux \cite{Capa},
Le\'on and Tudor \cite{LT},
 Nualart and Vives \cite {NV}, or Picard \cite{1pi}).
 Moreover, in this space, we can take
advantage of the pathwise characterizations of some stochastic integrals,
as we do in this paper, to deal with applications of the stochastic analysis
(see Le\'on et al. \cite{LSV}, Picard \cite{Pi} or Privault \cite{Pr}).
 In particular, the 
gradient operator  is a difference one.

Recently, several approaches to develope a calculus of variations for L\'evy 
processes have been introduced in some articles (see, for instance,
Di Nunno et al.  \cite{Nupro}, L{\o}kka \cite{Lo} and Sol\'e et al.
\cite{SUV}, among others).  The gradient and divergence operators are the 
fundamental tools in this theory again. In this paper, we restrict ourselves to
the canonical L\'evy space defined in \cite{SUV} because, in this space,
the  gradient operator defined utilizing the chaotic 
decomposition of a square--integrable random variable is not a ``derivative 
operator'' (see Section 2.3 below), but it is the sum of a derivative and an 
increment quotient
 operators. This fact is important because we can obtain and use
the relation between the stochastic integral introduced via the chaos 
decomposition and the pathwise stochastic integral, both with respect 
to the jump part of the involved L\'evy process (see Lemma \ref{L11.1}
below).

 The purpose of this paper is to use the Malliavin calculus on the
canonical L\'evy space given in \cite{SUV} to prove an anticipating It\^o
formula for L\'evy processes. Here,
the stochastic integrals with respect to the continuous and jump
parts of the underlying L\'evy process are in the Skorohod and pathwise sense,
respectively.
The coefficients in this formula
have two ``derivatives in the future sense''. It means, they are in a
class of square-integrable processes $u$ such that $u_t$ is in the domain
of the gradient operator $D$ at time $r$ for $r>t$, and $D_ru_t$ is also in
the domain of $D$ (see Section \ref{sub:2.4}). An example of processes
satisfying this property is the square-integrable and adapted processes,
whose ``derivative'' is equal to zero.

The paper is organized as follows. In Section \ref{sec:2} we present the 
framework that we use in this paper, Namely, we introduce some basic
facts of the canonical L\'evy space and of the Malliavin calculus on
this space. Finally, the anticipating It\^o formula is studied in Section
\ref{sec:3}.

\section{Preliminaries}

\label{sec:2} In this section we give the framework that will be used in
this article. That is, 
we introduce briefly  the It\^{o} multiple
integrals with respect to a L\'evy process, and
the canonical L\'{e}vy process 
considered by Sol\'{e}, Utzet and Vives \cite{SUV}.
 Then we present some basic facts on the Malliavin calculus for this process.
 We need to study the anihilation and creation
operators corresponding to the Fock space associated with the chaos
decomposition  on L\'{e}vy space, and
analyze the Sobolev spaces associated with these operators. Althoug some of
 these facts are known, we give them for the convenience of the reader. 

Throughout, we set $\erre_0 =\erre -\{0\}$ and $T>0$. 
Let $\nu$ be a L\'evy measure
 on $\erre$ such that $\nu (\{0\})=0$ and $\int_{\errita}x^2 d\nu (x) <\infty$
 (see Sato \cite{S}). The Borel $\sigma$--algebra of a set
 $A\subset {\hbox {{\gordas R}}}$ is
denoted by ${\cal B} (A)$. The jump of a c\'adl\'ag process $Z$ at time 
$t\in[0,T]$ is represented by $\Delta Z_t$ (i.e., $\Delta Z_t:=Z_t-Z_{t-}$).

\subsection{It\^{o} multiple integrals}

\label{sub:2.2} The construction of multiple integrals with respect to
L\'evy processes is quite similar to that of  multiple integrals with
respect to the Brownian motion. The reader can consult It\^o \cite{ito} for a
complete survey on this topic.

Let $X=\{X_{t}:t\in \lbrack 0,T]\}$ be a L\'{e}vy process with triplet $%
(\gamma ,\sigma ^{2},\nu )$. It is well--known that $X$ has the L\'{e}vy--It%
\^{o} representation (see \cite{S})
\begin{equation}\label{eq:truncada}
X_{t}=\gamma t+\sigma W_{t}+\int_{(0,t]\times
\{|x|>1\}}xdJ(s,x)
+\lim_{\varepsilon \downarrow 0}\int_{(0,t]\times
\{\varepsilon <|x|\leq 1\}}xd\widetilde{J}(s,x).
\end{equation}
Here the convergence is with probability 1, uniformly on $t\in[0,T]$,
 $W=\{W_{t}:t\in \lbrack 0,T]\}$ is a standard Brownian motion, 
\[
J(B)=\#\{t:(t,\Delta X_{t})\in B\},\;\;B\in {\cal B}([0,T]
\times {\hbox {{\gordas R}}}_{0}),
\]
is a Poisson measure with parameter $dt\otimes d\nu $ 
 and $d\widetilde{J}(t,x)
=dJ(t,x)- dtd\nu(x)$.

For $E_{1},\ldots ,E_{n}\in {\cal B}([0,T]\times {\hbox
{{\gordas R}}})$ such that $E_{i}\cap E_{j}=\emptyset $, $i\neq j$, and 
\[
\mu (E_{i}):=\sigma ^{2}\int_{\{t\in [0,T]:(t,0)\in
E_{i}\}}dt+ \int_{E_{i}-(E_i\cap ([0,T]\times\{0\}))}
 x^2 dtd\nu (x) <\infty,
\]
we define the multiple integral $I_{n}(1_{E_{1}\times \cdots \times E_{n}})$ 
of order $n$
with respect to $M$ by 
\begin{equation}
I_{n}(1_{E_{1}\times  \cdots \times  E_{n}})=M(E_{1})\cdots M(E_{n}),
\label{eq:4.3}
\end{equation}
with
$$
M(E_{i})=\sigma \int_{\{t\in [0,T]:(t,0)\in
E_{i}\}}dW_{t}+
\lim_{m\rightarrow \infty }\int_{\{(t,x)\in E_{i}:
\frac{1}{m}<|x|<m\}}xd%
\widetilde{J}(t,x), 
$$
where the limit is in the $L^2(\Omega)$ sense.
By linearity, we can define the multiple integral of order $n$
of an elementary function $f$ of the form
$$f(\cdot)=\sum_{i_1,\ldots,i_n=1}^N a_{i_1,\ldots,i_n}
1_{A_{i_1},\ldots,A_{i_n}}(\cdot),$$
where $A_1,\ldots,A_N$ are pairwise disjoint sets of 
${\cal B}([0,T]\times {\hbox{{\gordas R}}})$ and $a_{i_1,\ldots,i_n}=0$
if two of the indices $i_1,\ldots,i_n$ are equal.

The multiple integral $I_{n}$ is extended to $L_{n}^{2}:=L^{2}(([0,T]\times {
\hbox
{{\gordas R}}})^{n}$; ${\cal  B}(([0,T]\times {\hbox {{\gordas
R}}})^{n});\mu ^{\otimes n})$ due to the fact that 
the space of all the elementary functions
is dense in $L_{n}^{2}$ and  the property 
\begin{eqnarray}\label{eq:4.2}
\lefteqn{E[I_{n}(1_{E_{1}\times \cdots \times E_{n}})I_{m}(1_{F_{1}\times
\cdots \times F_{m}})]}  \nonumber   \\
&=&\delta _{n}(m)n!\int_{([0,T]\times \errita)^{n}}\widetilde{1}
_{E_{1}\times \cdots \times E_{n}}\widetilde{1}_{F_{1}\times \cdots \times
F_{m}}d\mu ^{\otimes n},
\end{eqnarray}
where $\widetilde{f}$ is the symmetrization of the function $f$ and
$\delta_n$ is the Dirac measure concentrated at $n$.

It is well--known that if $F$ is a square--integrable random variable,
measurable with respect to the filtration generated by $X$, then $F$ has the
unique representation 
\begin{equation}
F=\sum_{n=0}^{\infty }I_{n}(f_{n}),  \label{eq:4.1}
\end{equation}
where $I_{0}(f_{0})=f_{0}=E(F)$ and $f_{n}$ is a symmetric function in $
L_{n}^{2}$. This is the so called chaotic representation property
for L\'evy processes.

\subsection{Canonical L\'evy space}

\label{sub:2.1} The purpose of this subsection is to present some basic
elements of the structure of the canonical L\'{e}vy space
on the interval $[0,T]$. For a more
detailed account of this subject, we refer to Sol\'{e}, Utzet and Vives \cite
{SUV}.

The construction of the canonical L\'evy space is divided in three steps, as
follows: \vglue.25cm \noindent \textbf{Step 1.} Here we introduce the
canonical space for a compound Poisson process. Toward this end, let  
  $Q$ be a probability measure on ${\hbox {{\gordas R}}}$,
supported on $S\in {\cal B}({\hbox {{\gordas R}}}_{0})$, and
$\lambda>0$.
Set 
\[
\Omega _{T}=\bigcup_{n\geq 0}([0,T]\times S)^{n},
\]
with $([0,T]\times S)^{0}=\{\alpha \}$, where $\alpha $ is an arbitrary
point. The set $\Omega_T$ is equipped with the $\sigma$--algebra 
\[
\mathcal{F}_T=\{B\subset \Omega_T: B\cap ([0,T]\times S)^n \in 
{\mathcal B} (([0,T]\times S)^n), \quad \hbox{\rm for all}\quad n\geq
1\}. 
\]
The probability $P_T$ on $(\Omega_T, \mathcal{F}_T)$ is given by 
\[
P_T(B\cap ([0,T]\times S)^n)=e^{-\lambda T} \frac{\lambda^n (dt\otimes
Q)^{\otimes n}(B\cap ([0,T]\times S)^n)}{n!}, 
\]
with $(dt \otimes Q)^0=\delta_\alpha$. Here $\delta_{\alpha}$ is the 
Dirac measure concentrated at $\alpha$.

The space $(\Omega_T, \mathcal{F}_T, P_T)$ is called the canonical space for
the compound Poisson process with L\'evy measure $\lambda Q$. A similar 
definition for the Poisson process was given in 
 Neveu \cite{Ne}, and Nualart and Vives \cite{NV}. 
In $(\Omega_T, {\cal F}_T, P_T)$ the process 
\[
X_t (\omega) =\left\{ 
\begin{array}{lcl}
\sum\limits^n_{j=1} x_j 1_{[0,t]}(t_j), & \mathrm{if} & \omega =((t_1, x_1),
\cdots
, (t_n, x_n)), \\ 
0, & \mathrm{if} & \omega =\alpha,
\end{array}
\right. 
\]
is a compound Poisson process with intensity $\lambda$ and jump law given by
the probability measure $Q$.   
\vglue .25cm 
\noindent 
\textbf{Step 2.} Now we
consider the canonical space for a pure jump L\'evy process with 
L\'evy measure $\nu$.

Let $S_{1}=\{x\in {\hbox {{\gordas R}}}:\varepsilon _{1}<|x|\}$ and 
$S_{k}=\{x\in
 {\hbox {{\gordas R}}}:\varepsilon
_{k}<|x|\leq \varepsilon _{k-1}\}$ for $k>1$. Here
$\{\varepsilon _{k}:k\geq 1\}$ is a strictly decreasing sequence of
positive numbers such that $\varepsilon _{1}=1$,  $\lim\limits_{k
\rightarrow \infty}
 \varepsilon _{k}=0$ and  $\nu(S_k)\neq 0$.
 Note that the fact that $\nu$ is a L\'evy measure implies that 
$\nu (S_k)<\infty$ for every $k\ge 1$.
 Now, the canonical L\'{e}vy space with measure
 $\nu$ is defined as 
\[
(\Omega _{J},\mathcal{F}_{J},\mathcal{P}_{J})=\bigotimes _{k\geq 1}(\Omega
^{(k)},\mathcal{F}^{(k)},P^{(k)}),
\]
where $(\Omega ^{(k)},\mathcal{F}^{(k)},P^{(k)})$ is the canonical space for
the canonical compound Poisson process $\{ X_t^{(k)}: t\in[0,T]\}$
with intensity $\lambda _{k}=\nu (S_{k})$ and
probability measure $Q_{k}=\frac{\nu (\cdot \cap S_{k})}{\nu (S_{k})}$. In
this case, for $\omega =(\omega ^{k})_{k\geq 1}\in\Omega_J$
 and $t\in \lbrack 0,T]$,
the limit 
\[
J _{t}(\omega )=\lim_{n\rightarrow \infty
}\sum_{k=2}^{n}(X_{t}^{(k)}(\omega ^{k})-t\int_{S_{k}}xd\nu (x
))+X_{t}^{(1)}(\omega ^{1})
\]
exists with probability $1$ and it is a pure jump L\'{e}vy process
with L\'{e}vy measure $\nu $. 
\vglue.25cm 
\noindent 
\textbf{Step 3.} The
canonical L\'{e}vy space on $[0,T]$ with L\'{e}vy measure $\nu $ is 
\[
(\Omega ,\mathcal{F},P)=(\Omega _{W}\otimes \Omega _{J},\mathcal{F}
_{W}\otimes \mathcal{F}_{J},P_{W}\otimes P_{J}),
\]
where $(\Omega _{W},\mathcal{F}_{W},P_{W})$ is the canonical Wiener space.
Here, for $\omega =(\omega ^{\prime },\omega ^{^{\prime \prime }})\in \Omega
_{W}\otimes \Omega _{J}$, the process 
\begin{equation}
X_{t}(\omega )=\gamma t+\sigma \omega ^{\prime }(t)+J_{t}(\omega ^{^{\prime
\prime }})  \label{eq:3.1}
\end{equation}
is a L\'{e}vy process with triplet $(\gamma ,\sigma ^{2},\nu )$. For this
fact we refer to Sato \cite{S}.

\subsection{The anihilation and  creation operators}

\label{sub:2.3} Henceforth we suppose that the underlying probability space $%
(\Omega, \mathcal{F}, P)$ is the canonical L\'evy space with L\'evy measure $%
\nu$ and that $X$ is the L\'evy process defined in (\ref{eq:3.1}).

We say that the square-integrable
random variable $F$ given by (\ref{eq:4.1}) belongs to the
domain of the anihilation operator $D\;(F\in {\hbox {{\gordas D}}}^{1,2}$ for
short) if and only if
\begin{equation}
\sum_{n=1}^{\infty }nn!||f_{n}||_{L^2_n}^{2}<\infty .  \label{eq:5.1}
\end{equation}
In this case we define the random field $DF=\{D_zF:z\in \lbrack
0,T]\times {\hbox {{\gordas R}}}\}$ as 
\[
D_zF=\sum_{n=1}^{\infty }nI_{n-1}(f_{n}(z,\cdot )).
\]
Note that (\ref{eq:5.1}) yields that the last series converges in $
L^{2}(\Omega\times [0,T]\times {\hbox {{\gordas R}}}; 
 P\otimes\mu)$ by (\ref{eq:4.2}). Thus, in this case, we have
that $\sum_{n=0}^mI_n(f_n)$ and $\sum_{n=1}^{m}nI_{n-1}(f_{n}(z,\cdot ))$
converge to $F$ and $DF$ in $L^2(\Omega)$ and in
$L^{2}(\Omega\times[0,T]\times {\hbox {{\gordas R}}};P\otimes\mu)$ as 
$m\rightarrow\infty$,
respectively.
$D$ is a closed operator from $
L^{2}(\Omega )$ into $L^{2}(\Omega\times[0,T]\times {\hbox {{\gordas R}}};
P\otimes\mu)$, with dense domain. Similarly we
can define the iterated derivative $D_{z_{1},\ldots
,z_{n}}^{n}=D_{z_{1}},\cdots D_{z_{n}}$ and its domain 
${\hbox {{\gordas D}}}^{n,2}$.

The following result is due to Sol\'{e} et al. \cite{SUV} and it establishes
how we can figure out the random field $DF$ without 
using the chaos decomposition
(\ref{eq:4.1}). In order to state it, we need
the following:

Henceforth $W=\{W_t:t\in [0,T]\}$ is the canonical Wiener process and ${%
\hbox {{\gordas D}}}^{1,2}_W (L^2 (\Omega_J))$ denotes the family of $
L^2(\Omega_J, \mathcal{F}_J, P_J)$--valued random variables that are in the
domain of the derivative operator $D^W$ with respect to $W$. The reader
can consult Nualart \cite{Nul} for the basic definitions and
properties of this operator. The space 
 ${%
\hbox {{\gordas D}}}^{1,2}_W (L^2 (\Omega_J))$ is constructed as follows.
We say that a random variable $F$ is an $L^2(\Omega_J)$-valued smooth 
random variable if it has the form
$$F=f(W_{t_1},\ldots,W_{t_n})Z,$$
with $t_i\in[0,T]$, $f\in C_b^{\infty}({\hbox {{\gordas R}}}^n)$
(i.e., $f$ and all its partial derivatives are bounded), and
$Z\in L^2(\Omega_J, \mathcal{F}_J, P_J)$. 
The derivative of $F$ with respect to $W$, in the Malliavin
calculus sense, is defined as
$$D^W F=\sum_{i=1}^n\frac{\partial f}{\partial x_i}(W_{t_1},\ldots,W_{t_n})
Z1_{[0,t_i]}.$$ 
It is easy to see that $D^W$ is a closeable operator from
$L^2(\Omega_W;L^2 (\Omega_J))$ into 
$L^2(\Omega_W\times [0,T];L^2 (\Omega_J))$. Thus we can introduce the
space ${\hbox {{\gordas D}}}^{1,2}_W (L^2 (\Omega_J))$ as the completion
of the $L^2(\Omega_J)$-valued smooth 
random variables with respect to the seminorm
$$|| F||_{1,2,W}^2=E\left[|F|^2+|DF|_{L^2([0,T]}^2\right] .$$

For $\omega =(\omega',(\omega ^{k})_{k\geq 1})\in\Omega$, with
$\omega ^{k}=((t^k_1,x^k_1),\ldots,(t_{n_k}^k,x^k_{n_k}))$,
$F\in L^2(\Omega)$ and
$ z=(t,x)\in
(0,T]\times S_{k_0}$, for some positive integer $k_0$, 
we define 
\[
(\Psi_{t,x} F) (\omega) =\frac{F(\omega_z)-F(\omega )}{x},
\]
with $\omega_z =(\omega',(\omega ^{k}_z)_{k\geq 1})$ and 
$$
\omega ^{k}_z=\left\{
\begin{array}{cl}
((t,x),(t_1^{k_0},x_1^{k_0})\ldots,(t_{n_{k_0}}^{k_0},
x_{n_{k_0}}^{k_0})),& \hbox{if }
k=k_0,\\
\omega ^{k},
\phantom{holassssssssssssssslalalqa}&\hbox{otherwise}.
\end{array}\right.$$

\begin{lemma}
\label{L6.1} Let $F\in L^2 (\Omega)$ be a random variable such that:

\begin{itemize}
\item[i)]  $F\in {\hbox {{\gordas D}}}^{1,2}_W (L^2 (\Omega_J))$.
\item[ii)]  $\Psi F\in L^2 (\Omega \times [0,T] 
\times {\hbox {{\gordas R}}}_0;P \otimes \mu)$.
\end{itemize}
Then $F\in {\hbox {{\gordas D}}}^{1,2}$ and 
\[
D_{t,x} F=1_{\{0\}} (x) \sigma^{-1}D^W_t F+1_{{\hbox {{\gorditas R}}}%
_0} (x) \Psi_{t,x} F. 
\]
\end{lemma}
\textbf{Proof.} The proof of this result is an immediate consequence of 
\cite{SUV} (Propositions 3.5 and 5.5). \hfill $\rule{5pt}{5pt}$

Now we establish an auxiliary tool needed for our results.

\begin{lemma}
\label{L6.2} Let $F\in {\hbox {{\gordas D}}}^{1,2}$. Then there exists a
sequence $\{F_n :n\geq 1\}$ of the form 
\begin{equation}  \label{eq:6.1}
F_n =\sum^{N}_{i=1} H_{i,n}Z_{i,n}
\end{equation}
such that:
\begin{itemize}
\item[i)]  $H_{i,n}$ is a smooth functional in $L^2 (\Omega_W)$ and $Z_{i,n}
\in {\hbox {{\gordas D}}}^{2,2} \cap L^\infty (\Omega_J)$.
\item[ii)]  $F_n$ (resp. $DF_n$) converges to $F$ (resp. $DF$) in $%
L^2(\Omega)$ (resp. $L^2 (\Omega \times [0,T]\times {\hbox {{\gordas R}}};
P\otimes \mu))$ as $n\rightarrow \infty$.
\end{itemize}
\end{lemma}

\noindent\textbf{Remarks}
\begin{itemize}
\item[i)] Observe that $N$ in equality (\ref{eq:6.1}) is a positive
integer depending only on $n$.
\item[ii)] By \cite{SUV} (Proposition 5.4),
 $\Psi
Z_{i,n}\in L^2 (\Omega \times [0,T] \times {\hbox {{\gordas R}}}_0; P \otimes 
\mu)$.
\end{itemize}

\noindent
\textbf{Proof.} Note that  it is enough to show the result holds for a
multiple integral of the form (\ref{eq:4.3}). That is 
\[
F=M(E_1)\cdots M(E_k), 
\]
where $E_1,\cdots, E_k$ are pairwise disjoint borel subsets of $
[0,T]\times \erre$. Indeed, in this case, the result is also true for a
random variable $G$ with a finite chaos decomposition because, by
the definition of the multiple integrals, there exists a sequence
$\{ G_m:m\ge 1\}$ of linear combinations of multiple integrals of the form
(\ref{eq:4.3}) such that $G_m\rightarrow G$ in $L^2(\Omega)$ and
$DG_m\rightarrow DG$ in $L^2(\Omega\times[0,T]\times{\hbox {{\gordas R}}};
P\otimes \mu)$,
as $m\rightarrow\infty$. Therefore, (\ref{eq:5.1}) implies that the result is 
satisfied.

Let $\varphi \in C^{\infty }({\hbox {{\gordas R}}})$ be a function such that 
\[
\varphi (x)=\left\{ 
\begin{array}{ll}
1, & |x|\leq 1, \\ 
0, & |x|\geq 2.
\end{array}
\right. 
\]
Set $\rho _{n}(x)=x\varphi (\frac{x}{n})$ and 
\[
F_{n}=\prod_{i=1}^{k}\left(\rho _{n}\biggl(\int_{\{s:(s,0)\in E_{i}\}}\sigma
dW_{s}\biggr)
+\rho _{n}\left( \lim_{m\rightarrow \infty }\int_{\{(s,y)\in E_{i}:
\frac{1}{m}<|y|<m\}}yd\widetilde{J}(s,y)\right)\right) .
\]
Then,
\begin{eqnarray*}
\lefteqn{\Psi_{t,x} \biggl(\rho_n \biggl(\lim_{m\rightarrow \infty} 
\int_{\{(s,y)\in E_i :\frac{1}{m} <|y|<m\}} y d\widetilde{J} (s,y)\biggr)
\biggr)}\\
&=& \frac{1}{x} \biggl(\rho_n \biggl(x 1_{E_i}(t,x)+
\lim_{m\rightarrow \infty}
\int_{\{(s,y)\in E_i : \frac{1}{m} <|y|<m\}}
y d\widetilde{J} (s,y)\biggr)\\
&&- 
 \rho_n \biggl(
\lim_{m\rightarrow \infty}\int_{\{(s,y)\in E_i : \frac{1}{m} <|y|<m\}}
y d\widetilde{J} (s,y)\biggr)\biggr) ,
\end{eqnarray*}
and
\begin{eqnarray*}
\lefteqn{\Psi_{r,z}\Psi_{t,x} \biggl(\rho_n \biggl(\lim_{m\rightarrow \infty} 
\int_{\{(s,y)\in E_i :\frac{1}{m} <|y|<m\}} y d\widetilde{J} (s,y)\biggr)
\biggr)}\\
&=& \frac{1}{xz} \biggl(\rho_n \biggl(x 1_{E_i}(t,x)+z1_{E_i}(r,z)+
\lim_{m\rightarrow \infty}
\int_{\{(s,y)\in E_i : \frac{1}{m} <|y|<m\}}
y d\widetilde{J} (s,y)\biggr)\\
&&- \rho_n \biggl(x 1_{E_i}(t,x)+
\lim_{m\rightarrow \infty}
\int_{\{(s,y)\in E_i : \frac{1}{m} <|y|<m\}}
y d\widetilde{J} (s,y)\biggr)\\
&&-\rho_n \biggl(z1_{E_i}(r,z)+
\lim_{m\rightarrow \infty}\int_{\{(s,y)\in E_i : \frac{1}{m} <|y|<m\}}
y d\widetilde{J} (s,y)\biggr)\\
&&+\rho_n \biggl(
\lim_{m\rightarrow \infty}\int_{\{(s,y)\in E_i : \frac{1}{m} <|y|<m\}}
y d\widetilde{J} (s,y)\biggr)\biggr)
\end{eqnarray*}
Hence, $\rho_n \biggl(\lim_{m\rightarrow \infty} 
\int_{\{(s,y)\in E_i :\frac{1}{m} <|y|<m\}} y d\widetilde{J} (s,y)\biggr)$ is
in ${\hbox {{\gordas D}}}^{2,2}$ due to  
\cite{SUV} (Lemma 5.2) or Lemma \ref{L6.1}.

Now the result follows from the facts that 
 $F_{n}\rightarrow F$ in $L^{2}(\Omega )$ as $n\rightarrow \infty $,
$|\rho _{n}(x)|\leq |x|$ and that there is a
constant $C$ independent of $n$ such that $|\rho _{n}^{\prime }(x)|+|\rho
_{n}^{^{\prime \prime }}(x)|\leq C$.
\hfill $\cajita$

An immediate consequence of the last two lemmas is the following:

\begin{corolary}\label{ideD12}
Let $F$ be a random variable in $ L^2(\Omega)$. 
Then $F\in {\hbox {{\gordas D}}}^{1,2}$ if and
only if $F\in {\hbox {{\gordas D}}}^{1,2}_W (L^2 (\Omega_J))$
and  $\Psi F\in L^2 (\Omega \times [0,T] 
\times {\hbox {{\gordas R}}}_0;P \otimes \mu)$.\end{corolary}
\textbf{Proof.} The proof follows from Lemmas \ref{L6.1} and \ref{L6.2},
and from \cite{SUV} (Proposition 4.8).\hfill $\cajita$
 
We will also need the following result.

\begin{lemma}
\label{L8.1} Let $F\in {\hbox {{\gordas D}}}^{1,2}$ be a bounded random
variable.
 Then $(FG) \in {\hbox {{\gordas D}}}%
^{1,2}$ for every $G$ of the form (\ref{eq:6.1}).
\end{lemma}
\textbf{Proof.} We first observe that $FG \in {\hbox {{\gordas D}}}
^{1,2}_W(L^2(\Omega_J))$ due to Corollary \ref{ideD12}.

Finally, we have 
\[
\Psi_{t,x} (FG)=(\Psi_{t,x}F) G+F\Psi_{t,x} G+
(F(\omega_{(t,x)})-F)\Psi_{t,x}G. 
\]
Therefore $\Psi (FG)\in L^2 (\Omega \times [0,T]\times \erre_0; P\otimes \mu)$.
 Consequently 
the proof is complete by Lemma \ref{L6.1}. \hfill $\rule{5pt}{5pt}$

The creation operator $\delta $ is the adjoint of $D:{\hbox {{\gordas D}}}%
^{1,2}\subset L^{2}(\Omega )\rightarrow L^{2}
(\Omega \times \lbrack 0,T]\times {\hbox
{{\gordas R}}}; P\otimes \mu )$. It means, $u$ belongs to Dom $\delta $ 
if and only if 
$u\in L^{2}(\Omega \times \lbrack 0,T]\times {\hbox {{\gordas R}}}; P\otimes
\mu )$ is such that there exists a square--integrable random variable $%
\delta (u)$ satisfying the duality relation 
\begin{equation}
E\biggl[\int_{[0,T]\times \errita} u(z)(D_{z}F)d\mu (z)\biggr]
=E[\delta (u)F],
\quad \hbox{\rm for every}\quad F\in {\hbox {{\gordas D}}}^{1,2}.
\label{eq:8.1}
\end{equation}
It is not difficult to show  that this duality relation gives
 that if $u$ has
the chaos decomposition 
\[
u(z)=\sum_{n=0}^{\infty }I_{n}(u_{n}(z,\cdot )),\quad z\in \lbrack 0,T]\times {%
\hbox
{{\gordas R}}},
\]
where $u_{n}\in L_{n+1}^{2}$ is a symmetric function in the last $n$
variables, then $\delta (u)$ has the chaos decomposition (see\cite{Nul})
\[
\delta (u)=\sum_{n=0}^{\infty }I_{n+1}({\tilde u}_{n}).
\]
The creation operator of a process multiplied by a random variable can be
calculated via the following two results, which have been considered
by Di Nunno et al. \cite{Nupro} for pure jump L\'evy processes.

\begin{proposition}
\label{P9.1} Let $F$ be a random variable as in Lemma \ref{L8.1} and $u\in 
\mathrm{Dom}\; \delta $ such that 
\[
E\biggl[\int_{[0,T]\times {\hbox {{\gorditas R}}}%
}(u(t,x)(F+xD_{t,x}F))^{2}d\mu (t,x )\biggr]<\infty .
\]
Then $(t,x)\mapsto u(t,x)(F+xD_{t,x}F)$ belongs to Dom $\delta $ if and only
if 
\[
\biggl(
F\delta (u)-\int_{[0,T]\times {\hbox {{\gorditas R}}}}u(t,x)D_{t,x}Fd\mu
(t,x)\biggr)\in L^{2}(\Omega ).
\]
In this case 
\[
\delta (u(t,x)F+xu(t,x)D_{t,x}F)=F\delta (u)-\int_{[0,T]\times {
\hbox
{{\gorditas R}}}}u(t,x)D_{t,x}Fd\mu (t,x).
\]
\end{proposition}
\textbf{Proof.} Let $G$ be a random variable as in
the right-hand side of (\ref{eq:6.1}). Then 
 Lemma \ref{L8.1} and its proof give
\begin{eqnarray*}
\lefteqn{E[GF\delta (u)]} \\
&=&E\biggl[\int_{[0,T]\times {\hbox {{\gorditas R}}}}u(t,x)D_{t,x}(FG)d\mu
(t,x)\biggr] \\
&=&E\biggl[\sigma ^{2}\int_{0}^{T}u(t,0)D_{t,0}(FG)dt+\int_{[0,T]\times {%
\hbox {{\gorditas R}}}_{0}}u(t,x)D_{t,x}(FG)d\mu (t,x)\biggr] \\
&=&E\biggl[\sigma ^{2}\int_{0}^{T}u(t,0)(D_{t,0}F)Gdt+\sigma
^{2}\int_{0}^{T}u(t,0)FD_{t,0}Gdt\biggr] \\
&&+E\biggl[\int_{[0,T]\times {\hbox {{\gorditas R}}}
_{0}}u(t,x)((D_{t,x}F)G+FD_{t,x} G+x(D_{t,x}F)D_{t,x}G)d\mu (t,x)\biggr] \\
&=&E\biggl[G\int_{[0,T]\times {\hbox {{\gorditas R}}}}u(t,x)D_{t,x}Fd\mu
(t,x)\biggr] \\
&&+E\biggl[\int_{[0,T]\times {\hbox {{\gorditas R}}}}(u(t,x)F+u(t,x)x
D_{t,x}F)D_{t,x}Gd\mu (t,x)\biggr].
\end{eqnarray*}
Therefore the proof is complete by Lemma \ref{L6.2} and by the duality
relation (\ref{eq:8.1}). \hfill $\rule{5pt}{5pt}$

The following result is an immediate consequence of the proof of Proposition 
\ref{P9.1}.

\begin{corolary}
\label{C1.4} Let $u$ and $F$ be as in Proposition \ref{P9.1}. Moreover
assume that $(t,x)\mapsto u(t,x)xD_{t,x}F$ belongs to Dom $\delta $. Then $
Fu\in \Dom \;\delta $ if and only if 
\begin{equation}
F\delta (u)-\delta (u(t,x)xD_{t,x}F)-\int_{[0,T]\times {\hbox {{\gorditas
R}}}}u(t,x)D_{t,x}Fd\mu (t,x)  \label{eq:10.1}
\end{equation}
is a square--integrable random variable. In this case $\delta (Fu)$ is equal
to (\ref{eq:10.1}).
\end{corolary}

\subsection{Sobolev spaces}\label{sub:2.4}
In this subsection we proceed as in Al\`os and Nualart \cite{AN} in order
 to define the spaces that contain the integrands in our It\^o
formula. 

 Let $\mathcal{S}_{T}$ be the family of processes of the form 
$u(\cdot )=\sum_{j=1}^{n}F_{j}h_{j}(\cdot )$, where $F_{j}$ is a random
variable of the form (\ref{eq:6.1}) and $h:[0,T]\times {\hbox {{\gordas R}}}%
\rightarrow {\hbox {{\gordas R}}}$ is a bounded measurable function. Note
that the fact that $\int_{{\hbox {{\gorditas R}}}}x^{2}d\nu (x)<\infty $
implies that $h\in L^{2}([0,T]\times {\hbox {{\gordas R}}}; \mu )$. Denote by 
${\hbox {{\gordas L}}}^{1,2,f}$ the closure of $\mathcal{S}_{T}$ with
respect to the seminorm 
$$
||u||^2_{1,2,f}=E\int_{[0,T]\times \errita}u(z)^{2}d\mu (z) 
+E\int_{\Delta _{1}^{T}}(D_{s,y}u(t,x))^2d\mu(s,y)d\mu (t,x),
$$
where 
\[
\Delta _{1}^{T}=\biggl\{((s,y),(t,x))\in ([0,T]\times {\hbox {{\gordas R}}}%
)^{2}:s\geq t\biggr\}.
\]
A random field $u= \{u(s,y):(s,y)\in \lbrack 0,T]\times {\hbox {{\gordas
R}}}\}$ in ${\hbox {{\gordas L}}}^{1,2,f}$ belongs to the space ${\hbox
{{\gordas L}}}^{1,2,f}_-$ if there is $D^{-}u\in L^{2}(\Omega \times \lbrack
0,T]\times {\hbox
{{\gordas R}}};P\otimes \mu )$ such that 
\[
\lim_{n\rightarrow \infty }\int_{0}^{T}\int_{{\hbox {{\gorditas R}}}
}\sup_{(s-\frac{1}{n})\vee 0 \leq r< s,y\leq x\leq y+\frac{1}{n}
}E[|D_{s,y}u(r,x)-D^{-}u(s,y)|^{2}]d\mu (s,y)=0.
\]
The random field $D^{-}u$ has been introduced in \cite{Nul} for the Wiener case,
and in \cite{Nupro} for the pure jump case.

The next result will be a useful tool to state the It\^o formula for the
operator $\delta$. Remember that  we are using the notation 
$\Delta  X_s =X_s -X_{s-}$.

\begin{lemma}
\label{L11.1} Let $u=\{u(s,x):(s,x)\in[0,T]\times {\hbox
{{\gordas R}}}\}$ be a measurable random field and 
$\varepsilon_1>\varepsilon >0$ such that:
\begin{itemize}
\item[i)]  There exists a constant $c>0$ such that $|u(s,y)|<c$,
 for all $(s,y) \in [0,T] \times\{\varepsilon<|x|\le \varepsilon_1 \}$.
\item[ii)]  For any sequences $\{s_{n}\in [0,s):n\in {\hbox {{\gordas N}}}\}$ 
and $%
\{y_{n}\in\{\varepsilon<|x|\le\varepsilon_1 \}
:n\in {\hbox {{\gordas
N}}}\}$  that converge to $s\in[0,T]$ and 
$y\in \{\varepsilon<|x|\le \varepsilon_1 \}$, 
respectively, we have that the limit
\[
u(s-,y)=\lim_{n,m\rightarrow \infty }u(s_{n},y_{m})
\]
is well--defined.
\item[iii)]  $u(\cdot -,\cdot )\in {\hbox {{\gordas L}}}_{-}^{1,2,f}$.
\end{itemize}
Then 
\begin{eqnarray*}
\lefteqn{\sum_{0<s\leq t}u(s-,\Delta X_{s})\Delta X_{s}1_{\{\varepsilon
<|\Delta X_{s}|\leq \varepsilon_1\}}} \\
&\quad =&\delta ((u(s-,y)+yD^{-}u(s-,y))1_{\{\varepsilon <|y|\leq
\varepsilon_1\}}1_{[0,t]}(s))\\
&&+\int_{0}^{t}\int_{\{\varepsilon <|y|\leq \varepsilon_1\}}
u(s-,y)yd\nu (y)ds \\
&&+\int_{0}^{t}\int_{\{\varepsilon <|y|\leq \varepsilon_1\}}D^-u(s-,y)d\mu
(s,y),\quad t\in \lbrack 0,T].
\end{eqnarray*}
\end{lemma}
\textbf{Proof.} The definition of the space ${\hbox {{\gordas L}}}^{1,2,f}$
implies that there exists a sequence $\{u^{(m)}\in \mathcal{S}_{T}:m\in {%
\hbox
{{\gordas N}}}\}$ such that 
\begin{eqnarray}
\lefteqn{E\left[(u(t-,y)-u^{(m)}(t,y))^{2}\right.}\nonumber\\
&\quad&\left.+\int_{t}^{T}\int_{\errita}
(D_{s,x}(u(t-,y)-u^{(m)}(t,y)))^{2}d\mu (s,x)\right]\rightarrow 0,
\label{eq:12.1}
\end{eqnarray}
as $m\rightarrow \infty $, for $\mu$-a.a. $(t,y)\in \lbrack 0,T]\times {%
\hbox
{{\gordas R}}}$. Hence we can choose a sequence $\
{\mathcal{A}}
_{n}=\{(s_{i}^{(n)},y_{j}^{(n)}):i,j\in \{1,\ldots ,N\}\}$ such that:

\begin{itemize}
\item[$\bullet$] $N$ is a positive integer that depends on $n$ and
goes to $\infty$ as $n\rightarrow\infty$.
\item[$\bullet$]  $0\leq s^{(n)}_1< \cdots < s^{(n)}_N\leq T,\; 
-\varepsilon_1\leq
y^{(n)}_1<y^{(n)}_2<\cdots <y^{(n)}_N\leq \varepsilon_1$.
\item[$\bullet $]  $0=\lim\limits_{n\rightarrow \infty }s_{1}^{(n)},T=\lim
s_{N}^{(n)},-\varepsilon_1
=\lim\limits_{n\rightarrow \infty }y_{1}^{(n)}$ and $%
\varepsilon_1=\lim\limits_{n\rightarrow \infty }y_{N}^{(n)}\cdot$ 
\item[$\bullet$] $\max\limits_i
(s_{i+1}^{(n)}- s_{i}^{(n)})\rightarrow 0$ and $\max\limits_i
(y_{i+1}^{(n)}-y_{i}^{(n)})\rightarrow 0$ as $n\rightarrow \infty $.
\item[$\bullet$] Property  (\ref{eq:12.1}) holds when we write $(s^{(n)}_i,
y^{(n)}_{j+1})$ instead of $(t,y)$.
\end{itemize}

Thus, from the duality relation (\ref{eq:8.1}), Proposition 
\ref{P9.1}, 
(\ref{eq:12.1}) and \cite{SUV} (Theorem 6.1), we obtain 
\begin{eqnarray*}
\lefteqn{\sum_{i,j=1}^{N-1}
 u(s^{(n)}_i- , y_{j+1}^{(n)})
\int_{]s_i^{(n)},s_{i+1}^{(n)}]}
\int^{y^{(n)}_{j+1}}_{y^{(n)}_j} y 
1_{\{\varepsilon <|y|\leq \varepsilon_1\}}1_{[0,t]}(s)d\widetilde{J} (s,y)}\\
&\quad&=\sum_{i,j=1}^{N-1}
u( s^{(n)}_i -, y^{(n)}_{j+1}) \delta  \left(1_{\{\varepsilon <|y|\leq
\varepsilon_1\}}
1_{]s_{i}^{(n)},s_{i+1}^{(n)}]}(s)
1_{]y_{j}^{(n)},y_{j+1}^{(n)}]}(y)
1_{[0,t]}(s)\right)
\\
&\quad&=\sum_{i,j=1}^{N-1}
\biggl\{  \delta \left(
1_{[0,t]}(s)
1_{\{\varepsilon <|y|\leq \varepsilon_1\}}
1_{]s_{i}^{(n)},s_{i+1}^{(n)}]}(s)
1_{]y_{j}^{(n)},y_{j+1}^{(n)}]}(y)
\right.\\
&\quad&\left.\quad\quad\times
(u(s_{i}^{(n)}-,y_{j+1}^{(n)})+yD_{s,y}u(s_{i}^{(n)}-,y_{j+1}^{(n)})
)\right)
 \\
&\quad&\ \,+\int_{s_{i}^{(n)}}^{s_{i+1}^{(n)}}
\int_{y_{j}^{(n)}}^{y_{j+1}^{(n)}}
1_{[0,t]}(s)1_{\{\varepsilon <|y|\leq
\varepsilon_1\}}D_{s,y}u(s_{i}^{(n)}-,y_{j+1}^{(n)})d\mu (s,y)\biggr\} .
\end{eqnarray*}
Indeed, by Proposition \ref{P9.1} we have that the last equality holds
when we change $u(s_i^{(n)}-,y_{j+1}^{(n)})$ by
$u^{(m)}(s_i^{(n)},y_{j+1}^{(n)})$. Consequently, we prove that our claim is
 true using (\ref{eq:8.1})
with a random variable as in the right-hand side of (\ref{eq:6.1}) and 
letting $m$ go to $\infty$.
So, we can conclude the proof because of the dominated convergence theorem, the
hypotheses of this lemma and the fact that $\delta $ is a closed operator.
\hfill $\rule{5pt}{5pt}$

The space ${\hbox {{\gordas L}}}_{F}$ is the closure of $\mathcal{S}_{T}$
with respect to the norm 
\[
||u||_{F}^{2}=||u||_{1,2,f}^{2}+E\int_{\Delta
_{2}^{T}}(D_{r,x}D_{s,y}u(t,z))^{2}d\mu (r,x)d\mu (s,y) d\mu (t,z),
\]
with $\Delta _{2}^{T}=\{((r,x),(s,y),(t,z))\in ([0,T]\times {%
\hbox {{\gordas
R}}})^{3}:r\vee s\geq t\}$.

The following result was stated  on the Wiener space by Al\`os and Nualart 
\cite{AN}.

\begin{lemma}
\label{L4.1} Let $u\in {\hbox {{\gordas L}}}_F$. Then $u\in \mathrm{Dom}\;
\delta$ and 
\begin{equation}  \label{eq:14.1}
E[\delta (u)^2]\leq 2 ||u||^2_F.
\end{equation}
\end{lemma}
\textbf{Proof.} We first observe that it is enough to show that 
(\ref{eq:14.1}) is true for $u\in \mathcal{S}_T$ because $\delta$ is a closed
operator. In this case, we have by \cite{SUV} (Section 6) or by \cite{AN},
 
\begin{eqnarray}
E[\delta (u)^2] &=&E\left[\int_{0}^{T}\int_{{\hbox {{\gorditas R}}}
}u(t,x)^{2}d\mu (t,x)\right.  \nonumber  \label{eq:14.2} \\
&&+\left.\int_{0}^{T}\int_\errita\int_{0}^{T}
\int_{{\hbox{{\gorditas R}}}}D_{s,y}u(t,x)D_{t,x}u(s,y)d\mu (t,x)d\mu (s,y)\right].
\phantom{eeee}
\end{eqnarray}
Observe that 
\begin{eqnarray*}
\lefteqn{E\biggl[\int_{0}^{T}\int_{\errita}\int_{0}^{T}\int_{\errita}D_{s,y}u
(t,x)D_{t,x}u(s,y)d\mu (t,x)d\mu (s,y)\biggr]} \\
&=&2E\biggl[\int_{0}^{T}\int_{{\hbox {{\gorditas R}}}}u(s,y)\delta
( 1_{[0,s]}D_{s,y}u)d\mu (s,y)\biggr] \\
&\leq &E\biggl[\int_{0}^{T}\int_{{\hbox {{\gorditas R}}}}u(s,y)^2d\mu
(s,y)\biggr]+E\biggl[\int_{0}^{T}\int_{{\hbox {{\gorditas R}}}}[\delta
(1_{[0,s]} D_{s,y}u)]^{2}d\mu (s,y)\biggr] \\
&\leq &E\biggl[\int_{0}^{T}\int_{{\hbox {{\gorditas R}}}}u(s,y)^2d\mu
(s,y)\biggr]\\
&&+E\biggl[\int_{0}^{T}\int_{{\hbox {{\gorditas R}}}}\int_{0}^{s}
\int_{{\hbox
{{\gorditas R}}}}(D_{s,y}u(t,x))^{2}d\mu (t,x)d\mu (s,y)\biggr] \\
&&+E\biggl[\int_{0}^{T}\int_{{\hbox {{\gorditas R}}}}\int_{([0,s]\times {\hbox
{{\gorditas R}}})^{2}}D_{t,x}D_{s,y}u(r,z)d\mu (r,z)d\mu (t,x)d\mu 
(s,y)\biggr].\end{eqnarray*}
Thus (\ref{eq:14.2}) yields that (\ref{eq:14.1}) holds.
\hfill $\rule{5pt}{5pt}$

Inequality (\ref{eq:14.1}) allows us to consider Lemma \ref{L11.1} with
$\varepsilon=0$ or $\varepsilon_1=\infty$
 to obtain the relation between the pathwise integral
and the operator $\delta$.

\begin{corolary}\label{rel-total-sko-pathwise}
Let $u$ satisfy the hypotheses of Lemma \ref{L11.1} for each 
$\varepsilon,\varepsilon_1\in(a,b)$, with
$0\leq a$ and $b\le\infty$. Moreover assume that the random fields
$(s,y)\mapsto u(s-,y),yD^-
u(s-,y)$ belong to ${\hbox {{\gordas L}}}_F$ and $(s,y)\mapsto u(s-,y)y$
is pathwise integrable with respect to $\widetilde J$ on 
$[0,T]\times\{a<|y|<b\}$. Then
\begin{eqnarray*}
\lefteqn{\int_{]0,t]}\int_{\{a<|y|< b\}}u(s-,y)yd{\widetilde J}(s,y)}\\
&=&\delta\left((u(s-,y)+yD^-u(s-,y))1_{[0,t]}(s)1_{\{a<|y|< b\}}
(y)\right)\\
&&+\int_0^t\int_{\{a<|y|< b\}}D^-u(s-,y)d\mu(s,y),\quad t\in[0,T].
\end{eqnarray*}
\end{corolary}
\textbf{Proof.} The result is an immediate consequence of Lemmas
\ref{L11.1} and \ref{L4.1}.\hfill $\rule{5pt}{5pt}$

\section{The It\^o formula}

\label{sec:3} Here we assume that, for $i\in \{1,\ldots ,n\}$, 
\setcounter{equation}{0}
\begin{eqnarray*}
Y_{t}^{(i)} &=&Y_{0}^{(i)}+\int_{0}^{t}u_{i}(s)dW_{s}+\int_{0}^{t}\sigma
_{s}^{(i)}ds+\int_{]0,t]}\int_{\{|x|>1\}}v_{i1}(s-,x)xdJ(s,x) \\
&&+\int_{]0,t]}\int_{\{0<|x|\leq 1\}}v_{i2}(s-,x)xd\tilde{J}(s,x),\quad
t\in \lbrack 0,T].
\end{eqnarray*}
The stochastic integrals with respect to $W$ and $J$ are in the Skorohod and
pathwise sense, respectively, and
\begin{itemize}
\item[(H1)]  $Y^{(i)}_0\in {\hbox {{\gordas D}}}^{1,2}$.
\item[(H2)]  $u_i\in {\hbox {{\gordas L}}}_{F}$ 
is such that $\{\int_0^tu_i(s)dW_s:t\in[0,T]\}$ has continuous paths
and there is a constant $M>0$
such that $\int_{0}^{T}u_{i}(s)^{2}ds\leq M$ with probability $1$.
\item[(H3)]  $\sigma^{(i)} \in {\hbox {{\gordas L}}}^{1,2,f}$ and 
$\int_{0}^{T}(\sigma^{(i)}_{s})^{2}ds\leq M$ with probability 1, for
some positive constant $M$.
\item[(H4)]  $v_{i1}$ 
 satisfies the assumptions of Corollary \ref
{rel-total-sko-pathwise} for $a=1$ and $b=\infty$. Moreover assume that
there is a positive constant $M$ such that $|v_{i1}|<M$ for
$(s,y)\in[0,T]\times\{1<|x|<\infty\}$.
\item[(H5)]  The hypotheses of Corollary \ref{rel-total-sko-pathwise}
hold for
$v_{i2}$ with $a=0$ and $b=1$, and 
there is a positive constant $M$ such that $|v_{i2}(s-,y)|\le M$, for
$(s,y)\in[0,T]\times \{0\le|x|\le 1\}$. Moreover assume that $D^-v_{i2}\in
{\hbox {{\gordas L}}}^{1,2,f}$.
\end{itemize}

Observe that by Lemma \ref{L4.1} and Corollary 
\ref{rel-total-sko-pathwise}, we have that 
$$
\int_{]0,t]}\int_{\{0<|x|\leq 1\}}v_{i2}(s-,x )xd\widetilde{J}(s,x)
$$
belongs to $L^{2}(\Omega )$, for all $t\in[0,T]$.
 Also observe that in \cite{AN} (Theorem 1) we can
find sufficient conditions that guarantee the continuity of the
stochastic integral $\{\int_0^tu_i(s)dW_s:t\in[0,T]\}$.

To show our It\^o formula, we first need to assume that our L\'evy
process defined in (\ref{eq:3.1}) has no small side jumps. So,
for $\varepsilon >0$, we need to use the notation 
\begin{eqnarray}
Y_{t}^{(i),\varepsilon }
&=&Y_{0}^{(i)}+\int_{0}^{t}u_{i}(s)dW_{s}+\int_{0}^{t}\sigma
_{s}^{(i)}ds+\int_{]0,t]}\int_{\{|x|>1\}}v_{i1}(s-,x)xdJ(s,x)  \nonumber
\label{eq:16.1} \\
&&+\int_{]0,t]}\int_{\{\varepsilon <|x|\leq 1\}}v_{i2}(s-,x)xd\widetilde{J}%
(s,x),\quad t\in  [0,T].
\end{eqnarray}
The $i$--th jump time of the compound Poisson process 
$
\{\int_{]0,t]}\int_{\{\varepsilon <|x|\}}xdJ(s,x):t\in[0,T]\}
$
is denoted by $T_{i}^{\varepsilon }$. We also use the notation 
$T_{0}^{\varepsilon }=0$.

\begin{theorem}
\label{T16.1} Assume that (H1)--(H5) hold, for $i\in \{1,\ldots ,n\}$,
 and that $F\in C_{b}^{2}({%
\hbox
{{\gordas R}}}^{n})$. Then,  the processes 
\begin{eqnarray*}
\lefteqn{\left(
\partial_{i}F(Y_{s-})(u_{i}(s)1_{\{y=0\}}+v_{i2}(s-,y)1_{\{0<|y|\leq 1\}})
\right.} \\
&&\quad\left.
+y1_{\{0<|y|<1\}}D^{-}(v_{i2}\partial_{i}F(Y_{\cdot -}))(s,y)\right)1_{[0,t]}(s)
\end{eqnarray*}
belong to $\mathrm{Dom}\; \delta $ and 
\begin{eqnarray*}
\lefteqn{F(Y_{t})-F(Y_{0})} \\
&=&\delta \left(\left[\partial
_{i}F(Y_{s-})(u_{i}(s)1_{\{y=0\}}+v_{i2}(s-,y)1_{\{0<|y|\leq 1\}})\right. 
\right.\\
&&\left.\left.
 +y1_{\{0<|y|\leq 1\}}D^{-}(v_{i2}\partial _{i}F(Y_{\cdot -}))(s,y)\right]
1_{[0,t]}(s)\right) \\
&&+\frac{1}{2}\int_{0}^{t}\partial _{i}\partial
_{j}F(Y_{s})u_{i}(s)u_{j}(s)ds+\int_{0}^{t}\partial _{i}F(Y_{s})\sigma
_{s}^{(i)}ds \\
&&+\int_{0}^{t}\partial _{i}\partial
_{j}F(Y_{s})(D^{-}Y^{(j)})(s,0)u_{i}(s)ds \\
&&+\int_{0}^{t}\int_{\{0<|y|\leq 1\}}D^{-}(\partial
_{i}F(Y_{\cdot -})v_{i2})(s,y)d\mu (s,y) \\
&&+\sum_{0\leq s\leq t}\{F(Y_{s-}+\Delta Y_{s})-F(Y_{s-})-\partial
_{i}F(Y_{s-})v_{i2}(s-,\Delta X_{s})\Delta X_{s}\}1_{\{0< |\Delta
X_{s}|\leq 1\}} \\
&&+\sum_{0\leq s\leq t}(F(Y_{s-}+\Delta Y_{s})-F(Y_{s-}))1_{\{1<|\Delta
X_{s}|\}},\quad t\in \lbrack 0,T].
\end{eqnarray*}
Here we use the convention of summation over repeated indexes.
\end{theorem}

\noindent\textbf{Remark} By (\ref{eq:truncada}), we have
$\Delta Y_s1_{\{0<|\Delta X_s|\le 1\}}=v_{i2}(s-,\Delta X_{s})\Delta X_{s}
1_{\{0< |\Delta
X_{s}|\leq 1\}}$.
\vglue .25cm
\noindent 
\textbf{Proof.} We first observe that the process $Y^{(i),\varepsilon}$
given by (\ref{eq:16.1})  evolves as 
\[
Y^{(i),\varepsilon}_t=Y^{(i),\varepsilon}_{T_j^{\varepsilon}}
 +\int^t_{T_j^\varepsilon} 
u_i (s)
dW_s+\int^t_{T_j^\varepsilon} \sigma^{(i)}_s ds -
\int_{]T_j^{\varepsilon},t]}\int_{\{\varepsilon<|x|\leq 1\}} v_{i2}(s-, x)x 
\nu (dx)ds 
\]
 on the
stochastic interval $]T_j^\varepsilon, T_{j+1}^{\varepsilon}[$.
Consequently, proceeding as in \cite{AN} and
using that $W$ and $J$ are independent, and Corollary
\ref{rel-total-sko-pathwise},   we have that 
$1_{[0,t]}\partial_i F(Y) u_i$ belongs to $\mathrm{Dom}\; \delta^W$, for $i\in
 \{1,\ldots, n\}$ and 

\begin{eqnarray}  \label{eq:18.1}
F(Y^\varepsilon_t)-F(Y_0) &=&  
\sum_{i=1}^{\infty}\left(F(Y^{\varepsilon}_{t\wedge T^{\varepsilon}_{i}-})-
F(Y^{\varepsilon}_{t\wedge T^{\varepsilon}_{i-1}})\right)\nonumber\\
&&+\sum_{i=1}^{\infty}\left(F(Y^{\varepsilon}_{t\wedge T^{\varepsilon}_{i}})-
F(Y^{\varepsilon}_{t\wedge T^{\varepsilon}_{i}-})\right)\nonumber\\
&=&
\int^t_0 \partial_i F(Y^\varepsilon _s) u_i
(s) dW_s +\int^t_0 \partial _i F(Y^\varepsilon _s) \sigma^{(i)}_s ds 
\nonumber \\
&&- \int^t_0 \partial_i F(Y^\varepsilon _s) \int_{\{\varepsilon <|x|\leq
1\}} v_{i2} (s-,x) xd\nu (x) ds  \nonumber \\
&&+\frac{1}{2} \int^t_0 \partial_i \partial_j F(Y^\varepsilon _s) u_i (s)
u_j (s) ds  \nonumber \\
&&+ \int^t_0 \partial_i \partial_j F(Y^\varepsilon_s) (D^-
Y^{(j),\varepsilon}) (s,0) u_i (s)ds  \nonumber \\
&&+\sum_{0\leq s\leq t}(F(Y^\varepsilon_{s-}+\Delta Y^\varepsilon_s)
-F(Y^\varepsilon_{s-})),\quad t\in [0,T],
\end{eqnarray}
with 
\begin{eqnarray}  \label{eq:19.2}
D^- Y^{(j), \varepsilon} (s,0) &=& D_{s,0} Y_0^{(j)} 
+\int^s_0 D_{s,0} u_j (r)
dW_r+\int^s_0 D_{s,0} \sigma^{(j)}_r dr  \nonumber \\
&&+ \delta (D_{s,0} (v_{j2}(r-,y)+y D^-v_{j2} (r-,y)) 1_{\{\varepsilon
<|y|\leq 1\}}1_{[0,s]}(r))  \nonumber \\
&&+\delta (D_{s,0} (v_{j1} (r-,y)+y D^- v_{j1} (r-,y) ) 1_{\{1<|y|\}}
1_{[0,s]}(r))  \nonumber \\
&&+\int^s_0 \int_{\{\varepsilon <|y|\leq 1\}} D_{s,0} (D^- v_{j2}(r-,y))d\mu
(r,y)  \nonumber \\
&&+ \int^s_0 \int_{\{1<|y|\}} D_{s,0}(D^- v_{j1}(r-,y)) d\mu (r,y)\nonumber\\
&&+ \int^s_0 \int_{\{1<|y|\}}y D_{s,0} v_{j1} (r-,y) d\nu (y) dr.
\end{eqnarray}
Now we divide the proof in several steps. 
\vglue .25cm 
\noindent 
\textbf{Step 1.} Here we see that $Y^{(i),\varepsilon}_t \rightarrow 
Y_t^{(i)}$ in 
$L^2(\Omega)$ as $\varepsilon \downarrow 0$, for every $t\in [0,T]$.

It follows, from (\ref{eq:16.1}) and Lemma \ref{L11.1}, 
\begin{eqnarray}
\lefteqn{Y_{t}^{(i),\varepsilon }}\nonumber\\
&=&Y_{0}^{(i)}+\int_{0}^{t}u_{i}(s)dW_{s}+\int_{0}^{t}\sigma _{s}^{(i)}ds 
\nonumber  \label{eq:19.1} \\
&&+\delta \biggl((v_{i1}(s-,y)+yD^{-}v_{i1}(s-,y))
1_{\{1<|y|\}}1_{[0,t]}(s)\biggr)  \nonumber \\
&&+\delta \biggl((v_{i2}(s-,y)+yD^{-}v_{i2}(s-,y))
1_{\{\varepsilon<|y|\leq 1\}}1_{[0,t]}(s)\biggr)  \nonumber \\
&&+\int_{0}^{t}\int_{\{\varepsilon <|y|\leq 1\}}D^{-}v_{i2}(s-,y)d\mu
(s,y)+\int_{0}^{t}\int_{\{1<|y|\}}v_{i1}(s-,y)yd\nu (y)ds  \nonumber \\
&&+\int_{0}^{t}\int_{\{1<|y|\}}D^{-}v_{i1}(s-,y)d\mu (s,y).
\end{eqnarray}
Thus   our claim follows by Corollary \ref{rel-total-sko-pathwise}.
Indeed, by Lemma \ref{L4.1}, we have that
$$\delta \biggl((v_{i2}(s-,y)+yD^{-}v_{i2}(s-,y))
1_{\{\varepsilon<|y|\leq 1\}}1_{[0,t]}(s)\biggr) 
+\int_{0}^{t}\int_{\{\varepsilon <|y|\leq 1\}}D^{-}v_{i2}(s-,y)d\mu
(s,y)$$
converges in $L^2(\Omega)$ to the pathwise integral
$\int_0^t\int_{\{0<|y|\le 1\}} v_{i2}(s-,y)d{\tilde J}(s,y)$.
\vglue.25cm \noindent 
\textbf{Step 2.} Now we show
that $\partial _{i}F(Y^{\varepsilon }_{\cdot-})
v_{i2}(\cdot-,\cdot)$ is in ${\hbox {{\gordas L}}}%
^{1,2,f}_-$.

We first observe that (\ref{eq:14.1}), (\ref{eq:19.1}) 
and \cite{SUV} (Section 6) yield $Y^{(i),\varepsilon }\in {%
\hbox {{\gordas L}}}_{-}^{1,2,f}$, $i\in\{1,\ldots,n\}$.
Hence, $Y^{(i),\varepsilon }_{\cdot -}\in {%
\hbox {{\gordas L}}}_{-}^{1,2,f}$ due to $E[|Y^{(i),\varepsilon }_t-
Y^{(i),\varepsilon }_{t-}|]=0$, for $t\in[0,T]$, which follows
from (\ref{eq:16.1}). Thus, 
$D^-Y^{(i),\varepsilon }=D^-Y^{(i),\varepsilon }_{\cdot -}$.
 Therefore,
it is clear the fact that $F(Y^{\varepsilon})$ and $v_{i2}$ are bounded 
implies that 
\begin{eqnarray}\label{eq:19.1'}
D^{-}(\partial _{i}F(Y^{\varepsilon }_{\cdot -})v_{i2}(\cdot-,\cdot))(s,0) &=&
\partial _{i}\partial _{j}F(Y_{s-}^{\varepsilon})
v_{i2}(s-,0)D^{-}Y^{(j),\varepsilon }(s,0)  \nonumber \\
&&+\partial _{i}F(Y^{\varepsilon }_{s-})D^-v_{i2}(s-,0).
\end{eqnarray}
On the other hand, the definition of the operator
$\Psi$ leads to write, for $r>t$, 
\begin{eqnarray*}
\lefteqn{\Psi _{r,x}(\partial _{i}F(Y_{t-}^{\varepsilon })v_{i2}(t-,y))}
\\
&=&(\Psi _{r,x}\partial _{i}F(Y_{t-}^{\varepsilon}))
v_{i2}(t-,y)+\partial _{i}F(Y_{t-}^{\varepsilon })\Psi
_{r,x}v_{i2}(t-,y) \\
&&+ x(\Psi_{r,x}v_{i,2}(t-,y))\Psi_{r,x}\partial_i F(Y^\varepsilon_{t-})\\
&=&v_{i2}(t-,y)\frac{\partial _{i}F(Y_{t-}^{\varepsilon
}+xD_{r,x}Y_{t}^{\varepsilon })-\partial _{i}F(Y_{t-}^{\varepsilon })}{x}
+\partial _{i}F(Y_{t-}^{\varepsilon })D_{r,x}v_{i2}(t-,y) \\
&&+(\partial _{i}F(Y_{t-}^{\varepsilon }+xD_{r,x}Y_{t}^{\varepsilon
})-\partial _{i}F(Y^{\varepsilon }_{t-}))D_{r,x}v_{i2}(t-,y),
\end{eqnarray*}
which, together with (\ref{eq:19.1'})
and Corollary \ref{ideD12}, gives that 
$\partial _{i}F(Y^{\varepsilon })v_{i2}\in {
\hbox {{\gordas L}}}^{1,2,f}_-$, with
\begin{eqnarray*}
\lefteqn{D^-(\partial _{i}F(Y_{\cdot-}^{\varepsilon })v_{i2}(\cdot -,\cdot))(s,y)}
\\
&=&\left(\partial _{i}\partial _{j}F(Y_{s-}^{\varepsilon})
v_{i2}(s-,0)D^{-}Y^{(j),\varepsilon }(s,0)  
+\partial _{i}F(Y^{\varepsilon }_{s-})D^-v_{i2}(s-,0)\right)1_{\{y=0\}}\\
&&+\biggl(v_{i2}(s-,y)\frac{\partial _{i}F(Y_{s-}^{\varepsilon
}+yD^-Y^{\varepsilon }(s,y))-\partial _{i}F(Y_{s-}^{\varepsilon })}{y}
\\
&&+\partial _{i}F(Y_{s-}^{\varepsilon })D^-v_{i2}(s,y) \\
&&+(\partial _{i}F(Y_{s-}^{\varepsilon }+yD^-Y^{\varepsilon}(s,y)
)-\partial _{i}F(Y^{\varepsilon }_{s-}))D^-v_{i2}(s,y)\biggr)
1_{{\hbox {{\errita}}}_{0}}(y).
\end{eqnarray*}
\vglue.25cm 
\noindent 
\textbf{Step 3.} From
Step 2, Lemma \ref{L11.1} and (\ref{eq:18.1}), we get 
\begin{eqnarray}
\lefteqn{F(Y_{t}^{\varepsilon})}\nonumber\\
 &=&F(Y_{0})+\int_{0}^{t}\partial _{i}F(Y_{s}^{\varepsilon
})u_{i}(s)dW_{s}+\int_{0}^{t}\partial _{i}F(Y_{s}^{\varepsilon })\sigma
_{s}^{(i)}ds  \nonumber  \label{eq:21.1} \\
&&+\delta \left((\partial _{i}F(Y_{s-}^{\varepsilon
})v_{i2}(s-,y)\right.\nonumber\\
&&\quad \left.+y(D^{-}\partial _{i}F(Y^{\varepsilon
}_{\cdot -})v_{i2})(s,y))1_{\{\varepsilon <|y|\leq 1\}}1_{[0,t]}(s)\right)  
\nonumber \\
&&+\int_{0}^{t}\int_{\{\varepsilon <|y|\leq 1\}}D^{-}(\partial
_{i}F(Y_{\cdot -}^{\varepsilon })v_{i2})(s,y)d\mu (s,y)  \nonumber \\
&&+\frac{1}{2}\int_{0}^{t}\partial _{i}\partial _{j}F(Y_{s}^{\varepsilon
})u_{i}(s)u_{j}(s)ds  \nonumber \\
&&+\int_{0}^{t}\partial _{i}\partial _{j}F(Y_{s}^{\varepsilon
})(D^{-}Y^{(j),\varepsilon })(s,0)u_{i}(s)ds  \nonumber \\
&&+\sum_{0\leq s\leq t}(F(Y_{s-}^{\varepsilon }+\Delta Y_{s}^{\varepsilon
})-F(Y_{s-}^{\varepsilon })-\partial _{i}F(Y_{s-}^{\varepsilon
})\nonumber\\
&&\quad\quad\quad\times
v_{i2}(s-,\Delta X_{s})\Delta X_{s})1_{\{\varepsilon <|\Delta
X_{s}|\leq 1\}}  \nonumber \\
&&+\sum_{0\leq s\leq t}(F(Y_{s-}^{\varepsilon }+\Delta Y_{s}^{\varepsilon
})-F(Y_{s-}^{\varepsilon }))1_{\{1<|\Delta X_{s}|\}}.
\end{eqnarray}
\textbf{Step 4.} Now we analyze the convergence in $L^{2}(\Omega )$ of the
terms in (\ref{eq:21.1}). 
\begin{eqnarray*}
\lefteqn{E\biggl[|\sum_{0\leq s\leq t}(F(Y_{s-}^{\varepsilon }+\Delta
Y_{s}^{\varepsilon })-F(Y_{s-}^{\varepsilon }))1_{\{1<|\Delta
X_{s}|\}}|^{2}\biggr]} \\
&=&E\left[|\sum_{0\leq s\leq t}(F(Y_{s-}+\Delta Y_{s})-F(Y_{s-}
))1_{\{1<|\Delta X_{s}|\}}|^{2}\right] \\
&\leq &CE\left[\left(\sum_{i=1}^{n}\sum_{0\leq s\leq t}|v_{i1}(s-,\Delta
X_{s})\Delta X_{s}|1_{\{1<|\Delta X_{s}|\}}\right)^{2}\right] \\
&\leq &n^{2}CE\left[\left(\sum_{0\leq s\leq t}|\Delta X_{s}|1_{\{1<|\Delta
X_{s}|\}}\right)^{2}\right] \\
&\leq &CE\left[\left(\int_{]0,t]}\int_{\{|x|>1\}}|x| d\widetilde{J}
(s,x)+\int_{]0,t]}\int_{\{|x|>1\}}|x| d\nu (x)ds\right)^{2}\right] \\
&\leq &C\int_{]0,t]}\int_{\{|x|>1\}}x^{2}d\nu (x)ds+\biggl(
\int_{]0,t]}\int_{\{|x|>1\}}x d\nu (x)ds\biggr)^{2}\\
&\leq &C\int_{]0,t]}\int_{{\hbox {{\gorditas R}}}_{0}}x^{2}d\nu(x)ds<\infty
.
\end{eqnarray*}
Also
\begin{eqnarray*}
\lefteqn{E\biggl[\biggl(\sum_{0\leq s\leq t}(F(Y_{s-}+\Delta
Y_{s})-F(Y_{s-})-\partial _{i}F(Y_{s-})v_{i2}(s,\Delta X_{s})\Delta
X_{s})1_{\{0<|\Delta X_{s}|\leq \varepsilon \}}\biggr)^{2}\biggr]} \\
&\leq &E\biggl[\biggl(\sum_{i=1}^{n}\sum_{0\leq s\leq t}|\Delta
Y_{s}^{(i)}|^{2}1_{\{0<|\Delta X_{s}|\leq \varepsilon \}}\biggr)^{2}\biggr] \\
&=&E\biggl[\biggl(\sum_{i=1}^{n}\sum_{0\leq s\leq t}|v_{i2}(s-,\Delta
X_{s})\Delta X_{s}|^{2}1_{\{0<|\Delta X_{s}|\leq \varepsilon \}}\biggr)^{2}
\biggr]\\
&\leq &CE\biggl[\biggl(\sum_{0\leq s\leq t}|\Delta X_{s}|^{2}1_{\{0<|\Delta
X_{s}|\leq \varepsilon \}}\biggr)^{2}\biggr] \\
&\leq &CE\biggl[\biggl(\int_{]0,t]}\int_{\{0<|x|\leq \varepsilon \}}x^{2}d
\widetilde{J}(s,x)\biggl)^{2}\biggr]+C\biggl(
\int_{]0,t]}\int_{\{0<|x|\leq \varepsilon
\}}x^{2}d\nu (x)ds\biggr)^2 \\
&\leq &C\int_{]0,t]}\int_{\{0<|x|\leq \varepsilon \}}x^{2}d\nu
(x)ds\rightarrow 0\quad \mathrm{as}\quad \varepsilon \rightarrow 0.
\end{eqnarray*}
It is not difficult to  deduce, from Step 1, 
\[
E\left[\int_{0}^{t}|\partial _{i}F(Y_{s})u_{i}(s)-\partial
_{i}F(Y_{s}^{\varepsilon })u_{i}(s)|^{2}ds\right]\rightarrow 0
\]
and, from Step 2, (\ref{eq:14.1}) and the dominated convergence theorem, 
\begin{eqnarray*}
\lefteqn{E\biggl[\int_{]0,t]}\int_{\{0<|y|\leq 1\}}\biggl|\biggl(\partial
_{i}F(Y_{s}^{\varepsilon })v_{i2}(s-,y)+yD^{-}(\partial
_{i}F(Y^{\varepsilon })v_{i2})(s,y)\biggr)}\\
&\quad&\times 1_{\{\varepsilon <|y|\leq 1\}} 
-\partial _{i}F(Y_{s})v_{i2}(s-,y)+yD^{-}(\partial
_{i}F(Y)v_{i2})(s,y)\biggr|^{2}d\mu (s,y)\biggr] \\
&\quad&\quad\rightarrow 0\quad \mathrm{as}\quad \varepsilon \downarrow 0.
\end{eqnarray*}

The missing terms can be analyzed similarly.
\vglue.25cm
\noindent
\textbf{Step 5.} Finally the result follows from the fact that $\delta $ is a
closed operator and from Steps 1-4. \hfill $\rule{5pt}{5pt}$

\begin{theorem}
\label{T23.1} Assume that $\int_{{\hbox {{\gorditas R}}}_{0}}|x|d\nu
(x)<\infty $. Then  the hypotheses of Theorem \ref{T16.1} imply that 
\begin{eqnarray*}
F(Y_{t}) &=&F(Y_{0})+\int_{0}^{t}\partial
_{i}F(Y_{s})u_{i}(s)dW_{s}+\int_{0}^{t}\partial _{i}F(Y_{s})\sigma
_{s}^{(i)}ds \\
&&-\int_{0}^{t}\partial _{i}F(Y_{s})\int_{\{0<|x|\leq 1\}}v_{i2}(s-,x)xd\nu
(x)ds \\
&&+\frac{1}{2}\int_{0}^{t}\partial _{i}\partial _{j}F(Y_{s})u_{i}(s)u_j 
(s)ds \\
&&+\int_{0}^{t}\partial _{i}\partial
_{j}F(Y_{s})(D^{-}Y^{(j)})(s,0)u_{i}(s)ds \\
&&+\sum_{0\leq s\leq t}(F(Y_{s-}+\Delta Y_{s})-F(Y_{s-})),\quad t\in
\lbrack 0,T].
\end{eqnarray*}
\end{theorem}
\textbf{Proof.} The fact that $\int_{{\hbox {{\gorditas R}}}_{0}}|x|d\nu
(x)<\infty $ yields 
$$
E\biggl[\biggl(\int_{0}^{t}\int_{\{0<|x|\leq 1\}}|v_{i2}(s-,x)x|d\nu(x)ds%
\biggr)^{2}\biggr]
\leq C\biggr(\int_{0}^{t}\int_{\{0<|x|\leq 1\}}|x|d\nu (x)ds\biggr),
$$
which implies 
\begin{eqnarray*}
\lefteqn{E\biggl[\biggl(\int_{0}^{t}\partial _{i}F(Y_{s}^{\varepsilon
})\int_{\{\varepsilon <|x|\leq 1\}}v_{i2}(s-,x)x d\nu (x)ds} \\
&\quad&-\int_{0}^{t}\partial _{i}F(Y_{s})\int_{\{0<|x|\leq 1\}}v_{i2}(s-,x)x
d\nu(x)ds\biggr)^{2}\biggr]\rightarrow 0.
\end{eqnarray*}
Also we have 
\begin{eqnarray*}
\lefteqn{E\biggl[\biggl(\sum_{0\leq s\leq t}(F(Y_{s-}+\Delta
Y_{s})-F(Y_{s-}))1_{\{0<|\Delta X_{s}|\leq \varepsilon \}}\biggr)^{2}\biggr]} \\
&\leq &CE\biggl[\biggl(\sum_{i=1}^{n}\sum_{0\leq s\leq t}|v_{i2}(s-,\Delta
X_{s})\Delta X_{s}|1_{\{0<|\Delta X_s |\leq \varepsilon \}}\biggr)^{2}\biggr] \\
&\leq &CE\biggl[\biggl(\sum_{0\leq s\leq t}|\Delta X_{s}|1_{\{0<|\Delta X_{s}|\leq
\varepsilon \}}\biggr)^{2}\biggr] \\
&\leq &CE\biggl[
\biggl(\int_{]0,t]}\int_{\{0<|x|\leq \varepsilon \}}|x|d\widetilde{J}%
(s,x)+\int_{]0,t]}\int_{\{0|<|x|\leq \varepsilon \}}|x|d\nu (x)ds\biggr)^{2}
\biggr]\\
&\leq &C\int_{]0,t]}\int_{\{0<|x|\leq \varepsilon \}}x^{2}d\nu (x)ds+C\biggl(%
\int_{]0,t]}\int_{\{0<|x|\leq \varepsilon \}}|x|d\nu (x)ds\biggr)^{2}.
\end{eqnarray*}
Thus the result is a consequence of the proof of Theorem \ref{T16.1}. \hfill 
$\rule{5pt}{5pt}$
\vglue.25cm
\noindent\textbf{Acknowledgments.} Part of this paper was done while
Jorge A. Le\'on was visiting the Institut Mittag-Leffler. He is 
thankful for its hospitality.

\end{document}